\documentclass[11pt]{article}
\usepackage{palatino}

\usepackage{makeidx}
\usepackage{graphicx}

\usepackage{multicol}
\usepackage[centertags]{amsmath}
\usepackage{amsfonts}
\usepackage{euscript}
\usepackage{subfigure}
\usepackage{numbysec}
\usepackage{citesort}
\usepackage{color}
\allowdisplaybreaks[1]

\title{\textsf{Singular Behavior of Electric Field of \\ High Contrast Concentrated Composites}}
\author{Yuliya Gorb\thanks{%corresponding author,
Department of Mathematics, University of Houston, Houston, TX, 77204-3008, gorb@math.uh.edu}}
\date{}

\newtheorem{theorem}{Theorem}[section]
\newtheorem{lemma}[theorem]{Lemma}
\newtheorem{proposition}[theorem]{Proposition}

\begin{document}

\maketitle \thispagestyle{empty}

\begin{abstract}
\noindent A heterogeneous medium of constituents with vastly different mechanical properties, whose inhomogeneities are in close proximity to each other, is considered. The gradient of the solution to the corresponding problem exhibits singular behavior (blow up) with respect to the distance between inhomogeneities. This paper introduces a concise procedure for capturing the leading term of gradient's asymptotics precisely. This procedure is based on a thorough study of the system's energy. The developed methodology allows for straightforward generalization to heterogeneous media with a nonlinear constitutive description.
\end{abstract}

\section{Introduction}

This paper is on the study of blow up phenomena that occur in heterogeneous media consisting of a finite-conductivity matrix and perfectly conducting inhomogeneities (particles or fibers) close to touching. This investigation is motivated by the issue of material failure
initiation where one has to assess the magnitude of local fields,
including extreme electric or current fields, heat fluxes, and
mechanical loads, in the zones of high field concentrations. Such
zones are normally created by large gradient flows confined in very thin regions between particles of different potentials, see e.g. \cite{mark,bud-car,keller93,basl}.

%\textbf{Idea of this approach is behind thorough studies of the energy functional.}

These media are described by elliptic or degenerate elliptic equations with discontinuous coefficients. The problem of analytical study of solution regularity for such problems has been actively studied since 1999, and %these studies have 
resulted in series of papers \cite{bonnetier-vogelius,li-nirenbrg,li-vogelius,basl,akl,aklll,akllz,yun01,yun02,lim-yun,bly,bly2,gn12}
investigating different cases based on dimensions, shape of
inclusions, applied boundary conditions, etc. The main result up to
date can be summarized as follows: \textit{For two perfectly
conducting particles of an arbitrary smooth shape located at
distance $\delta$ from each other and away from the external
boundary, typically there exists $C>0$ independent of $\delta$
such that} 
\begin{equation}\label{E:prev}
\frac{1}{C\sqrt{\delta}}\leq \|\nabla u\|_{L^{\infty}}
\leq \frac{C}{\sqrt{\delta}} \quad \mbox{for} \quad d=2,~\qquad\frac{\log\delta^{-1}}{C~\delta}\leq \|\nabla u\|_{L^{\infty}}
\leq \frac{C \log\delta^{-1}}{\delta} \quad \mbox{for} \quad d=3,
\end{equation}
and corresponding bounds for the case of $N>2$ particles and $d>3$, see 
\cite{bly2}. %\cite{bly2} also provides similar to \eqref{E:prev} bounds for $d>3$, but the current paper focuses only on physically relevant dimensions  $d=2,3$. %, and some attempts to obtain similar estimates for insulating particles \cite{bly2,akl}. 
It is important to note that even though in some
referred studies it was mentioned on what parameters the constant
$C$ in \eqref{E:prev} depends upon, the precise asymptotics have not
been captured, only bounds for it have been established. Moreover,
methods in the aforementioned contributions have their limitations,
e.g. some of them use methods that work only in $2$D, some deal
with inhomogeneities of spherical shape only, and the
developed techniques, except one \cite{gn12} by the author, were
designed to treat \textit{linear} problems only, with no direct
extension or generalization to a nonlinear case. 

In the current paper an approach for gradient estimates for problems with particles of degenerate properties that
works for any number of particles of arbitrary shape in any dimensions is presented. The advantage and novelty is that the rate of blow up of  the electric field is captured {\bf precisely} as opposed to the existing methods and allows
for direct extensions to the nonlinear case (e.g. $p$-Laplacian). In particular, it is shown that
\begin{equation}\label{E:result}
 \|\nabla u\|_{L^{\infty}}=\frac{C}{\sqrt{\delta}} \quad \mbox{for} \quad d=2,~\qquad \|\nabla u\|_{L^{\infty}}
= \frac{C \log\delta^{-1}}{\delta} \quad \mbox{for} \quad d=3,
\end{equation}
with explicitly computable constant $C$ that depends on dimension $d$, particles array and their shapes, and an applied boundary field.

%\textbf{Idea of this approach is behind thorough studies of the energy functional.}

%\textcolor{red}{ The paper provides a very concise yet elegant procedure based on ... energy, that allows for {\it precise} capturing of the gradient's singular behavior}. 

The rest of the paper is organized as follows.  Chapter \ref{S:formulation} provides the problem setting and formulation of main results, proof of which is presented in Chapter \ref{S:proof}. Discussion of possible extensions is done in Chapter \ref{S:extensions} and conclusions are given in Chapter \ref{S:conclusions}. Proofs of auxiliary facts are shown in Appendices. 

\vspace{5pt}

\noindent \textbf{Acknowledgements}. The author thank A. Novikov for helpful discussions on the subject of the paper. 
%The author was supported by the NSF grant DMS-1016531.

\section{Problem Formulation and Main Results} \label{S:formulation}

%\subsection{Formulation} %\label{s:formulation}

\begin{figure}[ht]
 \centering
 \subfigure[Particles are at distance $\delta$ from each other]{
  \includegraphics[scale=0.8]{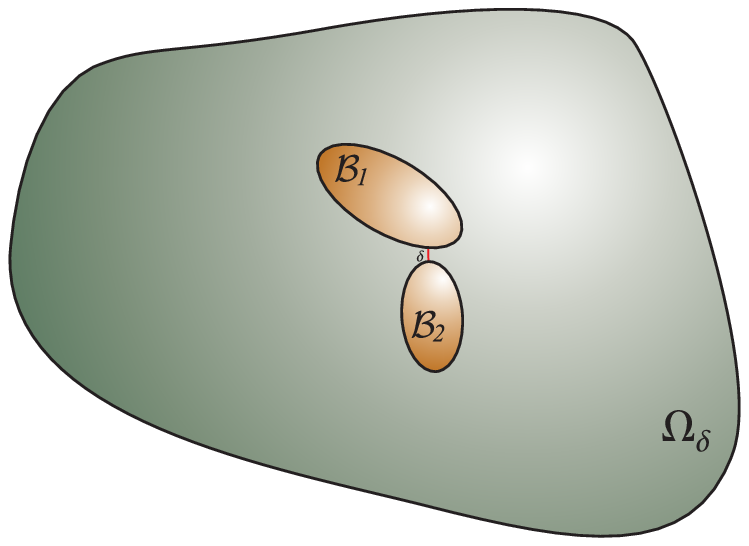}
   \label{F:domain2}
   }\qquad
 \subfigure[Particles are touching at one point]{
  \includegraphics[scale=0.8]{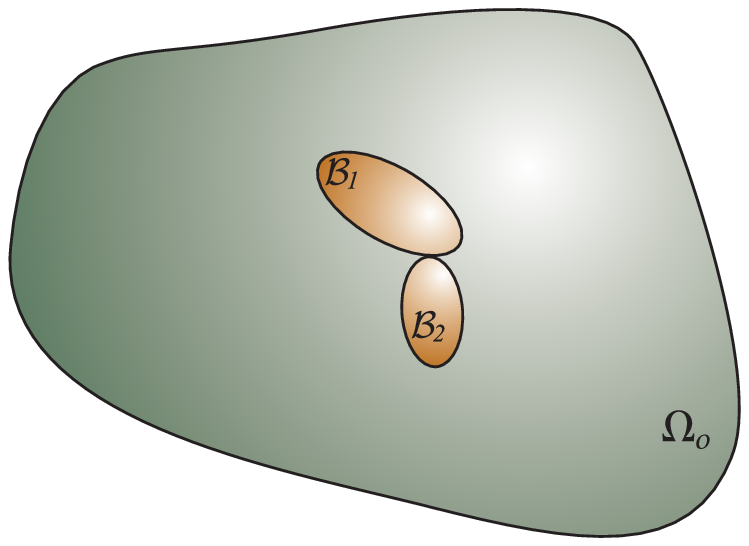}
  \label{F:domain0}
   }
\label{F:domain}
 \caption[Optional caption for list of figures]{Configurations of the composite occupying the domain $\Omega$ with particles $\mathcal{B}_1$ and $\mathcal{B}_2$}
\end{figure}

The current paper focuses only on physically relevant dimensions  $d=2,3$. To that end, let  $\Omega \in \mathbb{R}^d$, $d=2,3$ be an open bounded domain
with $C^{1,\alpha}$ ($0<\alpha\leq 1$) boundary
$\Gamma$.  It contains two particles
$\mathcal{B}_1$ and $\mathcal{B}_2$ %centered at $\boldsymbol{x}_i$, $i=1,\ldots,N$ of variable radii $R_i>0$
with smooth boundaries at the distances $0<\delta \ll 1$ from each
other; see Figure \ref{F:domain}.
We assume
\begin{equation}\label{E:distance}
\hbox{dist}(\partial \Omega, \mathcal{B}_1\cup \mathcal{B}_2) \geq K
\end{equation}
for some $K$ independent of $\delta$.  Let $\Omega_\delta$ model the matrix (or
the background medium) of the composite, that is,
$\Omega_\delta=\Omega/\overline{(\mathcal{B}_1\cup \mathcal{B}_2)}$, in which we
consider
\begin{equation} \label{E:pde-form-main}
\left\{
\begin{array}{r l l }
\triangle u(x) & = 0, & \displaystyle x\in \Omega_{\delta}\\[3pt]
 u(x) &=\mbox{const}, & \displaystyle  x\in \partial \mathcal{B}_i,~ i=1,2\\[3pt]
\displaystyle \int_{\partial \mathcal{B}_i}\frac{\partial u}{\partial n}  ~ds&=0, & i=1,2\\[3pt]
u(x)  &= U(x),  & \displaystyle  x\in  \Gamma
\end{array}
\right.
\end{equation}
where a bounded weak solution $u$ represents the electric
potential in $\Omega_\delta$, and $U$ is the given
applied potential on the external boundary $\Gamma$. Note
$u$ takes a constant value, that we denote $\mathcal{T}_i$, on the boundary of particle $\mathcal{B}_i$ ($i=1,2$). This is a unique constant for which the zero-flux condition, that is the third equation of \eqref{E:pde-form-main}, is satisfied. The constants $\mathcal{T}_1$, $\mathcal{T}_2$ are unknown apriori and should be found in the course of solving the problem.

The goal is to derive the asymptotics of the solution gradient with respect to the small parameter $\delta\ll 1$ that defines the close proximity of particles to each other. To formulate the main result of the paper, consider an auxiliary problem defined as follows. Construct a line connecting the centers of mass of  $\mathcal{B}_1$ and $\mathcal{B}_2$ and ``move'' particles toward each other along this line until they touch. Denote now  the newly obtained domain outside of particles by $\Omega_o$ where we consider the following problem:
\begin{equation} \label{E:pde-form-v0}
\left\{
\begin{array}{r l l }
\triangle v_o(x) & = 0, & x\in\Omega_o\\[3pt]
v_o(x) & =\mbox{const}, & \displaystyle x\in\partial \mathcal{B}_1\cup \partial\mathcal{B}_2\\[3pt]
\displaystyle \int_{\partial\mathcal{B}_1}\frac{\partial v_o}{\partial n}~ds+ \int_{\partial\mathcal{B}_2}\frac{\partial
v_o}{\partial n}~ds  & =0, \\[5pt]
v_o(x) & = U(x) ,  & x\in\Gamma
\end{array}
\right.
\end{equation}
This problem differs from \eqref{E:pde-form-main} by that the potential takes the {\it same} constant value on the boundaries of {\it both} particles. We denote this potential by $\mathcal{T}_o$ and introduce a number that depends on the external potential $U$:
\begin{equation} \label{E:R0}
\mathcal{R}_o=\mathcal{R}_o[U]:=%\mathcal{R}_0[U]:=
\int_{\partial\mathcal{B}_1}\frac{\partial v_o}{\partial n}~ds.
\end{equation}
Without loss of generality, we assume that particles potential in \eqref{E:pde-form-main} satisfy $\mathcal{T}_2>\mathcal{T}_1$, which would mean that $\mathcal{R}_o>0$ for sufficiently small $\delta$.
 
The following theorem summarizes the main result of this study. 
 
\begin{theorem} \label{T:main}
The asymptotics of the electric field for problem \eqref{E:pde-form-main} is given by
\begin{equation} \label{E:main-res}
\|\nabla u\|_{L^\infty(\Omega_\delta)}=\left[1+o(1)\right]
\begin{cases}
\displaystyle \frac{\mathcal{R}_o}{\mathcal{C}_{12}}\frac{1}{\delta^{1/2}}, & d=2\\[8pt]
\displaystyle \frac{\mathcal{R}_o}{\mathcal{C}_{12}}\frac{1}{ \delta |\ln \delta|}, & d=3
\end{cases}
\quad \quad \mbox{for }~\delta\ll 1,
\end{equation}
with $\mathcal{R}_o$ defined above in \eqref{E:R0} and explicitly computable constant $\mathcal{C}_{12}$ that depends on curvatures of particle boundaries $\partial \mathcal{B}_1$ and $\partial \mathcal{B}_2$  at the point of the closest distance %between $\mathcal{B}_1$ and $\mathcal{B}_2$ 
and defined below in \eqref{E:C0}.
\end{theorem}

\section{Proof of Main Results} \label{S:proof}

The proof of Theorem \ref{T:main} consists of ingredients collected in the following facts.

In \cite{gn12} using the method of barriers it was shown that the electric field of the system associated with the problem:
\begin{equation} \label{E:phi-eq}
\left\{
\begin{array}{r l l }
\triangle \phi(x) & = 0, & x\in\Omega_ \delta\\[3pt]
\phi(x) & =T_i, & \displaystyle x\in\partial \mathcal{B}_i,~i=1,2\\[3pt]
\phi(x) & = U(x) ,  & x\in\Gamma
\end{array}
\right.
\end{equation}
stated on the same domain $\Omega_\delta$ %and particles $\mathcal{B}_1$ and $\mathcal{B}_2$, and 
with the same boundary potential $U$ as in  the above problem \eqref{E:pde-form-main} is given by
\[%begin{equation} \label{E:phi-grad}
\|\nabla \phi \|_{L^\infty(\Omega_\delta)}=\frac{|T_2-T_1|}{\delta}[1+o(1)], \quad \mbox{for }~ \delta\ll 1.
\]%end{equation}
 In contrast to \eqref{E:pde-form-main}, the constants $T_1$ and $T_2$ in \eqref{E:phi-eq} are arbitrary, which implies the solution of  \eqref{E:phi-eq} may not satisfy the integral identities the flux of $u$ on $\partial \mathcal{B}_i$ as in \eqref{E:pde-form-main}. In particular, one has 

\begin{lemma} \label{L:elfield}
The asymptotics of the electric field of \eqref{E:pde-form-main} is as follows:
\begin{equation} \label{E:u-grad}
\|\nabla u \|_{L^\infty(\Omega_\delta)}=\frac{\mathcal{T}_2-\mathcal{T}_1}{\delta}[1+o(1)], \quad \mbox{for }~ \delta\ll 1,
\end{equation}
where $\mathcal{T}_2$ and $\mathcal{T}_1$ are the potentials on $\mathcal{B}_1$ and $\mathcal{B}_2$, respectively, for which the zero integral flux condition as in \eqref{E:pde-form-main} satisfied.
\end{lemma}
With \eqref{E:u-grad} the problem is reduced to finding the asymptotics of the potential difference $\mathcal{T}_2-\mathcal{T}_1$ in terms of the distance parameter $\delta$, given in the proposition. 

\begin{proposition} \label{L:potdif}
The asymptotics of the potential difference $\mathcal{T}_2-\mathcal{T}_1$ is given by:
\begin{equation} \label{E:potdif}
\mathcal{T}_2-\mathcal{T}_1=\frac{\mathcal{R}_o}{g_\delta}[1+o(1)], \quad \mbox{for }~ \delta\ll 1,
\end{equation}
where $\mathcal{R}_o$ is defined by \eqref{E:R0} and $g_\delta$ by:
\begin{equation} \label{E:g}
g_\delta=\begin{cases}
\displaystyle \mathcal{C}_{12}\delta^{-1/2}, & d=2\\[3pt]
\displaystyle \mathcal{C}_{12} |\ln \delta|, & d=3
\end{cases}
\end{equation}
with constant $\mathcal{C}_{12}$ introduced below in \eqref{E:C0} that depends on curvatures of particles boundaries at the point of their closest distance.
\end{proposition}

\noindent \textit{Proof of Proposition  \ref{L:potdif}.} \\
The method of proof is based on observation that asymptotics \eqref{E:potdif} of $\mathcal{T}_2-\mathcal{T}_1$ can be derived by investigating the energy associated with the system \eqref{E:pde-form-main} and defined by:
\begin{equation} \label{E:enE}
\mathcal{E}=\int_{\Omega_\delta}|\nabla u|^2~dx,
\end{equation}
where $u$ solves \eqref{E:pde-form-main}. %In particular, the following result holds.
A remarkable feature of problem \eqref{E:pde-form-main} is that 
potentials $\mathcal{T}_1$ and $\mathcal{T}_2$ are minimizers of the energy quadratic form of the potentials: % $T_1$ and $T_2$ of the problem \eqref{E:phi-eq}
\begin{equation} \label{E:IML}
\mathcal{E}=\min_{\{T_1,T_2\}} E(T_1,T_2), \quad  E(T_1,T_2)= \int_{\Omega_\delta}|\nabla \phi|^2~dx, \quad 
\mbox{where }~ \phi~~  \mbox{solves }~ \eqref{E:phi-eq}.
\end{equation}
This observation is the essence of the so-called Iterative Minimization Lemma, first introduced in \cite{bgn2}. 
Therefore, if we find an approximation of $\mathcal{E}$ for sufficiently small $\delta$, we would be able to derive an asymptotics for $\mathcal{T}_2-\mathcal{T}_1$ then. For the energy $\mathcal{E}$ the following holds true.

%Namely, the following facts are true regarding $\mathcal{E}$. 

%\begin{lemma} \label{L:itmin}
%The system's energy $\mathcal{E}$ can be written as follows:
%\[\mathcal{E}=\min_{\{T_1,T_2\}} E(T_1,T_2) = \min_{\{T_1,T_2\}} \int_{\Omega_\delta}|\nabla \phi|^2~dx,\] 
%where $\phi $ solves \eqref{E:phi-eq}.
%\end{lemma}

\begin{lemma} \label{L:energy}
The energy $\mathcal{E}$ of  \eqref{E:pde-form-main} can be written as
\begin{equation} \label{E:min-pr}
\mathcal{E}=\min_{\{t_1,t_2\}} \left[a_1t_1^2+a_2t_2^2+2b_1t_1+2b_2t_2+2c_{12}t_1t_2+C\right],
\end{equation}
with asymptotics of coefficients of the quadratic form $E$:
%\[
%\begin{array}{r r l }
%a_1=&a_2=&\mathcal{G}_\delta[1+o(1)], \\[3pt]
%b_1=&-b_2=&\mathcal{R}_o,\\[3pt]
%&c_{12}=&-\mathcal{G}_\delta[1+o(1)],
%\end{array}
%\]
\begin{equation} \label{E:coef-as}
a_1=a_2=g_\delta[1+o(1)], \quad
b_1=-b_2=\mathcal{R}_o[1+o(1)], \quad
c_{12}=-g_\delta[1+o(1)], \quad \mbox{for }~ \delta\ll 1,
\end{equation}
and  $\mathcal{R}_o$ given by \eqref{E:R0}, and $g_\delta$ by \eqref{E:g}.
\end{lemma}

This lemma is proven in Appendix \ref{A:en-lem}. 
Now substituting asymptotics  \eqref{E:coef-as} of coefficients to  \eqref{E:min-pr} and dropping the low order terms we define the quadratic form:
\[
\hat{E}(t_1,t_2)=g_\delta(t_2-t_1)^2-2\mathcal{R}_o(t_2-t_1),
\] 
whose minimizer $(\hat{t}_1,\hat{t}_2)$ provides asymptotics of the sought potential difference, namely,
\[
\mathcal{T}_2-\mathcal{T}_1=|\hat{t}_2-\hat{t}_1|[1+o(1)]=\frac{\mathcal{R}_o}{g_\delta}[1+o(1)],  \quad \mbox{for }~ \delta\ll 1.
\]
This concludes the proof of Proposition   \ref{L:potdif}. $\Box$

\vspace{7pt}

\noindent \textit{Proof of Theorem  \ref{T:main}.}\\ 
Asymptotic relations \eqref{E:u-grad}-\eqref{E:potdif} and definition \eqref{E:g} yield main result \eqref{E:main-res} for sufficiently small $\delta$. 
$\Box$

\section{Extensions} \label{S:extensions}

\begin{figure}[!ht]
  \centering
  \includegraphics[scale=.59]{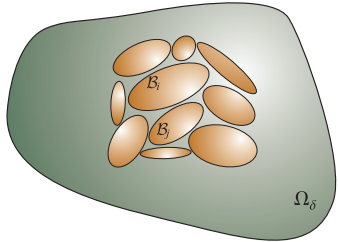}
  \caption{Composite containing $N>2$ particles $\mathcal{B}_1,\ldots,\mathcal{B}_N$}   \label{F:domainN}
\end{figure}

%\begin{remark}{}
\noindent {\bf 4.1. Extension to the case of $N>2$ particles}. \hspace{2pt}
The presented above approach allows for an extension to any number of particles $N>2$, where neighbors $\mathcal{B}_i$ and $\mathcal{B}_j$  are located at the distance $\delta_{ij}=O(\delta)\ll 1$ from each other, see Figure \ref{F:domainN}. The notion of ``neighbors'' can be defined based on% Delauney triangulation or 
the Voronoi tesselation with respect to the particles centers of mass, namely, the neighbors are the nodes that share the same Vonoroi face. % $x_i$, $i=1,\ldots,N$. 
In this case, similarly to above, one has to consider a ``limiting problem'' 
\eqref{E:pde-form-v0} in the domain $\Omega_o$ where the third condition is replaced to 
\[
\sum_{i=1}^N\int_{\partial\mathcal{B}_i}\frac{\partial v_o}{\partial n}~ds=0.\]
To obtain $\Omega_o$ one can connect centers of mass of neighboring pairs $\mathcal{B}_i$, $\mathcal{B}_j$ with lines, and ``move'' all particles alone those lines toward each other until $\partial \mathcal{B}_i$ touches at least one of its neighbor, where $i\in \{1,\ldots,N\}$, $j\in \mathcal{N}_i$, where $\mathcal{N}_i$ is the set of indices of neighbors to particle $\mathcal{B}_i$. 
Now, similarly to \eqref{E:R0}, introduce numbers
\[
\mathcal{R}_i=\mathcal{R}_i[U]:=
\int_{\partial\mathcal{B}_i}\frac{\partial v_o}{\partial n}~ds, \quad i\in \{1,\ldots,N\}.
\]
Then minimize the energy quadratic form $E$ as in \eqref{E:min-pr} and derive asymptotics of its coefficients in terms of $\mathcal{R}_i$ and $\delta_{ij}$ using $|\mathcal{T}_i-\mathcal{T}_\delta|\ll 1$ to obtain the potential difference asymptotics for the neighbors:
\begin{equation} \label{E:potdif-N}
|\mathcal{T}_i-\mathcal{T}_j|=\frac{|\mathcal{R}_i-\mathcal{R}_j|}{g_{ij}}[1+o(1)], \quad \mbox{for }~ \delta\ll 1, \quad \mbox{and }~ i\in \{1,\ldots,N\}, ~j\in \mathcal{N}_i.
\end{equation} 
Asymptotics of parameters $g_{ij}$ in \eqref{E:potdif-N} is similar to one of $g_\delta$ and given by
\[
g_{ij}=\mathcal{C}_{ij} \begin{cases}
\delta^{-1/2}_{ij},& d=2\\[3pt]
\left|\ln \delta_{ij}\right|,& d=3
\end{cases}, 
%\quad  \mathcal{C}_{ij}= \begin{cases} 2\pi \frac{\alpha_i\alpha_j}{\alpha_i+\alpha_j},& d=2\\[3pt] 4\pi \frac{a_ia_j}{a_i+a_j} \frac{b_ib_j}{b_i+b_j},& d=3 \end{cases}
\] 
with $\mathcal{C}_{ij}$ given by formula  \eqref{E:C0}  in Appendix \ref{A:coeff} %for $\mathcal{C}_{12}$ 
where $i$ should be replaced by $1$ and $j$ by $2$. 
Finally, use 
\[
\|\nabla u \|_{L^\infty(\Omega_\delta)}=\max_{i\in\{1,..,N\}, ~j\in \mathcal{N}_i}\frac{|\mathcal{T}_i-\mathcal{T}_j|}{\delta_{ij}}[1+o(1)], \quad \mbox{for }~ \delta\ll 1,
\]
with asymtotics \eqref{E:potdif-N} to obtain the blow up of electric field of the composite with more than two particles.
%\end{remark}

 \vspace{5pt}

%\begin{remark}{}
\noindent {\bf 4.2. Extension to the  nonlinear case}. \hspace{2pt}
One can also generalize the proposed methodology for high-contrast materials with the matrix described by nonlinear constitutive laws such as $p$-Laplacian. The system's energy in this case is given by $\displaystyle \mathcal{E}=\int_{\Omega_\delta}|\nabla u|^p~dx$, ($p>2$), where $u$ solves \eqref{E:pde-form-main} with first and third equations replaced by $\nabla \cdot \left( |\nabla u|^{p-2}\nabla u\right)=0$ in $\Omega_\delta$ and $\displaystyle \int_{\partial \mathcal{B}_i}|\nabla u|^{p-2} (\nabla u\cdot n) ~ds=0$, respectively. Note that for  a successful application of the described approach, one needs to show that the energy function $E(T_1,T_2)$, whose minimal value $ \mathcal{E}$ is attained at the solution $u$, is differentiable with respect to the potential $T_i$ on $\partial \mathcal{B}_i$.
The blow up of the electric field is then
\[
\|\nabla u\|_{L^\infty(\Omega_\delta)}=\left( \frac{\mathcal{R}_o}{\mathcal{C}_{12}} \right)^{\frac{1}{p-1}} \delta^{-\frac{d-1}{2(p-1)}}\left[1+o(1)\right], \quad p>2, \quad d=2,3, \quad \mbox{for }~\delta\ll 1,
\]
see also \cite{gn12}.
%\end{remark}

 \vspace{5pt}

%\begin{remark}{\bf (Extension to dimensions $d>3$)}
\noindent {\bf 4.3. Extension to dimensions $d>3$}. \hspace{2pt}
The described above procedure remains the same if one needs to obtain asymptotics for $|\nabla u|$ in dimensions greater than three. For this, one has to derive asymptotics of $g_\delta$ for $d>3$ first, following method described in Appendices \ref{A:G-lem} and \ref{A:coeff}. For simplicity of presentation we omit this case here.
%\end{remark}

\section{Conclusions} \label{S:conclusions}

%The goal of the paper is to establish precise estimates for ...

As observed in \cite{mark,bud-car,keller93}, in a composite consisting of a matrix of finite conductivity with perfectly conducting particles close to touching the electric field exhibits blow up. This blow up is, in fact, the main cause for a material failure which occurs in the thin gaps between neighboring particles of different potentials. The electric field of such composites is described by the gradient of the solution to the corresponding boundary value problem. The current paper provides a concise and elegant procedure for capturing the singular behavior of the solution gradient  {\it precisely} that does not require employing a heavy analytical machinery developed in previous studies  \cite{akl,aklll,akllz,yun01,yun02,lim-yun,bly,bly2,gn12}. This procedure relies on simple observations about energy of the corresponding system and its minimizers that were sufficient to acquire the sought asymtpotics exactly. 
%Moreover, the %tools %(
The techniques developed and adapted here are independent of dimension $d$, particles shape and their total number $N$,
whereas strict dependence on $d$ and particles shape was the main limitation of previous contributions on the subject  \cite{akl,aklll,akllz,yun01,yun02,lim-yun,bly,bly2,gn12}. Furthermore, the developed above procedure allows for a straightforward generalization to a {\it nonlinear} case. 
%In contrast to  \textcolor{red}{ \cite{} }, the generalization to ... of our approach does not require developing new ... tools. Last but not least, 

\section{Appendices}

\subsection{Proof of  Lemma \ref{L:energy}}   \label{A:en-lem}
%Here we prove Lemma \ref{L:energy}.\\
\textit{Proof.} Consider a family of auxiliary problems defined on the same domain $\Omega_\delta$ as \eqref{E:pde-form-main}:
\begin{equation} \label{E:pde-form-v-delta}
\left\{
\begin{array}{r l l }
\triangle v(x) & = 0, & \displaystyle x\in\Omega_\delta\\[3pt]
v(x) & =\mbox{const}, & \displaystyle x\in\partial \mathcal{B}_1\cup \partial\mathcal{B}_2\\[3pt]
\displaystyle \int_{\partial\mathcal{B}_1}\frac{\partial v}{\partial n}~ds  + \int_{\partial\mathcal{B}_2}\frac{\partial v}{\partial n}~ds& =0, \\[5pt]
v(x) & = U(x),  & \displaystyle x\in \Gamma
\end{array}
\right.
\end{equation}
%In contrast to the original \eqref{E:pde-form-main} 
As in  \eqref{E:pde-form-v0} the constant value of the potential is the same on both particles that we denote by $\mathcal{T}_\delta$. However, in contrast to \eqref{E:pde-form-v0} here particles are located at distance $\delta$ from each other while in \eqref{E:pde-form-v0} particles touch at one point. With that, similarly to \eqref{E:pde-form-v0} we introduce the number
\begin{equation} \label{E:R_delta}
\mathcal{R}_\delta[U]:=\int_{\partial\mathcal{B}_1}\frac{\partial v}{\partial n}~ds.
\end{equation}
In \cite{gn12}, it was shown that asymptotics of $\mathcal{R}_\delta[U]$ is given by
\begin{equation} \label{E:Rd-asym}
\mathcal{R}_\delta[U]=\mathcal{R}_o[1+o(1)], \qquad \mbox{for }~ \delta\ll 1.
\end{equation}
Using the linearity of problem  \eqref{E:pde-form-main} we decompose its solution into
\begin{equation} \label{E:u-decomp}
u=v+(\mathcal{T}_1-\mathcal{T}_\delta)\psi_1+(\mathcal{T}_2-\mathcal{T}_\delta)\psi_2,
\end{equation}
with $\psi_i$ ($i=1,2$) solving
\begin{equation} \label{E:psi-eq}
\left\{
\begin{array}{r l l }
\triangle \psi_i(x) & = 0, & x\in\Omega_ \delta\\[3pt]
\psi_i(x) & =\delta_{ij}, & \displaystyle x\in\partial \mathcal{B}_j,~~i,j\in\{1,2\}\\[3pt]
\psi_i(x) & = 0,  & x\in\Gamma
\end{array}
\right.
\end{equation}
where $\delta_{ij}$ is the Kroneker delta.
Invoking \eqref{E:u-decomp}, we compute the energy \eqref{E:enE} of the system and obtain:
\begin{equation} \label{E:en-decomp}
\mathcal{E}=\mathcal{E}_v+\mathcal{G}_1(\mathcal{T}_1-\mathcal{T}_\delta)^2+\mathcal{G}_2(\mathcal{T}_2-\mathcal{T}_\delta)^2+
2\mathcal{R}_\delta[U](\mathcal{T}_1-\mathcal{T}_2)+
2\mathcal{C}_{12}(\mathcal{T}_1-\mathcal{T}_\delta)](\mathcal{T}_2-\mathcal{T}_\delta),
\end{equation}
where 
\begin{equation} \label{E:not}
\mathcal{E}_v:=  \int_{\Omega_\delta}|\nabla v|^2~dx, \qquad 
\mathcal{G}_i:= \int_{\Omega_\delta}|\nabla \psi_i|^2~dx=\int_{\partial \mathcal{B}_i}\frac{\partial \psi_i}{\partial n}~ds, \quad i=1,2,
\end{equation}
are the energies of systems given by \eqref{E:pde-form-v-delta} and \eqref{E:psi-eq}, respectively, and 
\[
\mathcal{C}_{12}:= \int_{\Omega_\delta}(\nabla \psi_1 \cdot\nabla \psi_2)~dx.
\]
Trivial integration by parts yields that 
\begin{equation} \label{E:C12}
\mathcal{C}_{12}=-\mathcal{G}_1+C_1=-\mathcal{G}_2+C_2,%-\int_{\partial \mathcal{B}_1}\frac{\partial \psi_2}{\partial n}~ds=-\int_{\partial \mathcal{B}_2}\frac{\partial \psi_1}{\partial n}~ds.
\end{equation}
where constants $C_i$ depend on $d$, $K$ and shape of the particles, but independent of $\delta$. On the other hand, the problem \eqref{E:pde-form-v-delta} is regular in the sense that its electric field $|\nabla v|$ does not exhibit blow up since there is no potential drop between the particles. Hence, 
\begin{equation} \label{E:Ev}
\mathcal{E}_v=:C,
\end{equation}
 that depends on the same parameters as the above constants.
Finally,  in Appendix \ref{A:G-lem} we show that for sufficiently small $\delta$:
\begin{equation} \label{E:G-blowup}
\mathcal{G}_i=g_\delta[1+o(1)].
\end{equation}
With notations introduced in  \eqref{E:coef-as}, \eqref{E:Ev}, \eqref{E:not}, Iterative Minimization Lemma \eqref{E:IML}, and asymptotics \eqref{E:Rd-asym}, \eqref{E:C12}, \eqref{E:G-blowup}  we have  from \eqref{E:en-decomp}:
\[
\mathcal{E}=C+a_1(\mathcal{T}_1-\mathcal{T}_\delta)^2+a_2(\mathcal{T}_2-\mathcal{T}_\delta)^2+
2b_1\mathcal{T}_1+2b_2\mathcal{T}_2+
2c_{12}(\mathcal{T}_1-\mathcal{T}_\delta)(\mathcal{T}_2-\mathcal{T}_\delta),
\]
which with \eqref{E:IML} yields \eqref{E:min-pr}, where $t_i=T_i-\mathcal{T}_\delta$, $i=1,2$.\\
$\Box$

\subsection{Asymptotics of $\mathcal{G}_i$} \label{A:G-lem}

Here we prove asymptotic formula \eqref{E:G-blowup} which is stated in the following lemma.
\begin{lemma} \label{L:G-blowup}
For $\delta\ll 1$ asymptotics of  the energy $\mathcal{G}_i$ defined in \eqref{E:not} is given by
%\begin{equation} \label{E:G-blowup}
$\mathcal{G}_i=g_\delta[1+o(1)]$, $i=1,2$,
%\end{equation}
with $g_\delta$ defined in \eqref{E:g}.
\end{lemma}
\noindent \textit{Proof.} 
To derive an asymptotics of $\mathcal{G}_i$ we adopt the method of {\it variational bounds} that has become a classical tool in capturing the leading terms of asymptotics of the energy of the corresponding system. This method is based on two equivalent variational formulations of the corresponding problem that provide upper and lower bounds for the energy matching up the leading order of asymptotics. Employing this method we use of a couple observations made in \cite{bk,bn,bgn1} which are vital in capturing the sought asymtptotics. But before, we need to introduce a coordinate system in which the construction will be made.

\begin{figure}[ht]
 \centering
 \subfigure[]{
  \includegraphics[scale=1]{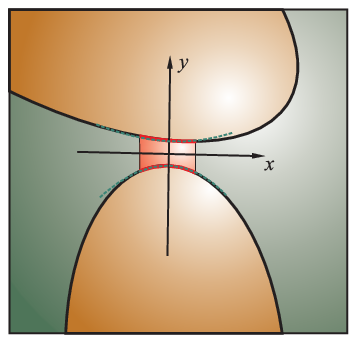}
   \label{F:neck01}
   }\qquad \qquad \qquad
 \subfigure[]{
  \includegraphics[scale=1]{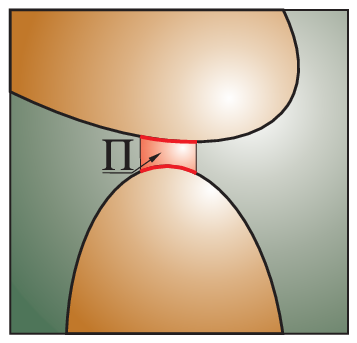}
  \label{F:neck02}
   }
\label{F:neck}
 \caption[]{(a) Coordinate system; \hspace{2pt} (b) Neck $\Pi$ between particles $\mathcal{B}_1$ and $\mathcal{B}_2$ }
\end{figure}

First, we write each point $ x\in  \mathbb{R}^d$ as $x=(\bar{x},x_d)$ where 
 \[\left\{
 \begin{array}{l l l}
\bar{x}=x,&x_d=y, & \mbox{when }~d=2\\[5pt]
\bar{x}=(x,y),&x_d=z, &  \mbox{when }~d=3
\end{array} \right.
 \]
Then, connect the centers of mass of particles with a line and ``move'' $\mathcal{B}_1$ and $\mathcal{B}_2$ along this line toward each other until they touch, thus, producing domain $\Omega_o$ as above in \eqref{E:pde-form-v0}. The point of their touching defines the origin of our cylindrical coordinate  system. The line connecting the centers will be the axis $Ox_d$, %$d=2,3$ (that is, $Ox_d=Oy$ when $d=2$, and $Ox_d=Oz$ when $d=3$), 
see Figure \ref{F:neck01}.
When particles are ``moved back'' at the distance $\delta$ from each other along $Ox_d$, we construct a ``cylinder'' of radius $w\gg \delta$ that contains this line. This ``cylinder'' is depicted as the red region in Figure \ref{F:neck02} that we call a {\it neck} and denote by  $\Pi$. Also, introduce the distance $H=H(\bar{x})$ between boundaries of $\mathcal{B}_1$  and $\mathcal{B}_2$, which in the selected coordinate system is a function of $\bar{x}\in  \mathbb{R}^{d-1}$.

%Now decompose domain $\Omega_\delta$ into the neck $\Pi$  and the rest of the domain. 
The mentioned above observations about energy estimates are as follows. First, the minimal value of the energy functional in the neck $\Pi$ is attained on the system with insulating lateral boundary $\partial \Pi$ of the cylinder, that is, 
\[
\mathcal{G}_i=\int_{\Omega_\delta}|\nabla \psi_i|^2~dx\geq \int_{\Pi}|\nabla \psi_i|^2~dx \geq \int_{\Pi}|\nabla \psi^i_{\Pi}|^2~dx,
\]
where function $ \psi^i_{\Pi}$ solves the problem %in the neck $\Pi$ with Neumann boundary conditions on the lateral boundary  $\partial \Pi$ of the cylinder:
\begin{equation} \label{E:psi-Pi}
\left\{
\begin{array}{r l l }
\triangle \psi^i_{\Pi}(x) & = 0, & x\in\Pi\\[3pt]
\psi^i_{\Pi}(x) & =\delta_{ij}, & \displaystyle x\in\partial \mathcal{B}_j,~~i,j\in\{1,2\}\\[5pt]
\displaystyle \frac{\partial \psi^i_{\Pi}}{\partial n}(x) & = 0,  & x\in \partial\Pi
\end{array}
\right.
\end{equation}
On the other hand, since energy $\mathcal{G}_i$ is the minimal value of the energy functional attained at the minimizer $\psi_i$, its upper bound is given by any test function $\phi_i$ from the set 
\[
V_i=\left\{\phi_i\in H^1(\Omega):~ \phi_i =\delta_{ij}~~\mbox{on } \partial \mathcal{B}_j,~\phi_i =0~\mbox{on }\Gamma\right\},
\]
via
\[
\mathcal{G}_i=\int_{\Omega_\delta}|\nabla \psi_i|^2~dx \leq \int_{\Omega_\delta}|\nabla \phi_i|^2~dx.
\]
Hence, the variational bounds for $\mathcal{G}_i$ are
\begin{equation} \label{E:en-bounds}
 \int_{\Pi}|\nabla \psi^i_{\Pi}|^2~dx\leq \mathcal{G}_i  \leq \int_{\Omega_\delta}|\nabla \phi_i|^2~dx, \quad \mbox{where } \phi_i \in V_i, ~\mbox{ and } \psi^i_{\Pi}\mbox{ solves } \eqref{E:psi-Pi}.
\end{equation}
Therefore, the problem is now reduced to construction of an approximation to $\nabla \psi^i_{\Pi}$ and finding a function $\phi_i\in V_i$ so that the integrals in \eqref{E:en-bounds} match up to the leading order for $\delta\ll 1$. For this purpose, one can use the {\it Keller's} functions \cite{keller63} defined in $\Pi$ by 
\begin{equation} \label{E:keller}
 \phi_\Pi^i(x)= \frac{x_d}{H(\bar{x})}, \quad x\in \Pi.
\end{equation}
With this $\phi_\Pi^i$, we define a test function $\phi_i\in V_i$ %that provides the upper variational bound 
by 
\[
\phi_i(x)=\begin{cases}
\phi_\Pi^i(x), & x\in \Pi\\[3pt]
\phi_o^i(x) &  x\in \Omega_\delta\setminus \Pi
\end{cases},
\]
where $\phi_o^i$ solves 
\[
\left\{
\begin{array}{r l l }
\triangle\phi_o^i(x) & = 0, & x\in \Omega_\delta\setminus \Pi\\[3pt]
\phi_o^i(x) & =\delta_{ij}, & \displaystyle x\in\partial \mathcal{B}_j,~~i,j\in\{1,2\}\\[5pt]
\phi_o^i(x) & = \phi_\Pi^i,  & x\in \partial\Pi\\[3pt]
\phi_o^i(x) & = 0,  & x\in \Gamma
\end{array}
\right.
\]
Employing the method of barriers to this problem one can show that $|\nabla \phi_o^i|\leq C$ with constant $C$ depending on $d$ and $K$ but independent of $\delta$. Thus,
\[
\mathcal{G}_i  \leq \int_{\Pi}|\nabla \phi_\Pi^i|^2~dx+C,\quad  \mbox{for } \delta\ll 1.
\]

The {\it dual variational principle} will help to estimate integral $\displaystyle \int_{\Pi}|\nabla \psi^i_{\Pi}|^2~dx$, namely, 
\[
\begin{array}{r l l }
\nabla \psi^i_{\Pi}& \displaystyle =\mbox{argmax}_{W_\Pi}\left[-  \int_{\Pi}j_i^2~dx+2\int_{\partial \mathcal{B}_i}(j_i\cdot n)~ds\right],\\[12pt]
W_{\Pi}& \displaystyle= \left\{j\in L^2(\Pi;\mathbb{R}^d):~\nabla\cdot j=0~\mbox{in }\Pi,~j\cdot n=0~\mbox{on } \partial\Pi\right\}.
\end{array}
\]
The test flux $j_i\in \mathbb{R}^d$ is chosen
\begin{equation} \label{E:test-flux}
j_i=\begin{cases}
\displaystyle \left(0,\frac{1}{H(x)}\right), & d=2\\[10pt]
\displaystyle  \left(0,0,\frac{1}{H(x,y)}\right), & d=3
\end{cases} 
\end{equation}
Therefore,
\begin{equation} \label{E:en-test-flux}
  \int_{\Pi}j_i^2~dx=\int_{\Pi}\frac{d\bar{x}}{H^2(\bar{x})}.
\end{equation}
Hence, we have two-sided bounds for $ \mathcal{G}_i$:
\[
- \int_{\Pi}j_i^2~dx+2\int_{\partial \mathcal{B}_i}(j_i\cdot n)~ds\leq \mathcal{G}_i\leq  \int_{\Pi}|\nabla \phi_\Pi^i|^2~dx +C,\quad  \mbox{for } \delta\ll 1.
\]
%We select such a function $\phi_i$ that the corresponding energy functional will coincide up to the leading order of asymptotics in $\delta$ with ... that provides the lower bound for $\mathcal{G}_i$. 
%If the distance between the two ...
%The test function in the neck $\Pi$ is 
%\[
%%%\phi_\Pi=\frac{x_d}{h(r)}, \quad x_d=\begin{cases}y, & d=2\\[3pt]z, & d=3\end{cases}, 
%%%\qquad r=\begin{cases} x, & d=2\\[3pt] \sqrt{x^2+y^2}, & d=3 \end{cases}.
%\phi_\Pi^i=\begin{cases}
%\displaystyle \frac{y}{H(x)}, & d=2\\[10pt]
%\displaystyle \frac{z}{H(x,y)}, & d=3
%\end{cases} 
With selected test functions $ \phi_\Pi^i$ and $j_i$ by \eqref{E:keller} and \eqref{E:test-flux}, respectively, it is trivial to show that the difference between the upper and lower bounds is simply % right- and left-hand sides of this estimate is %simply $\displaystyle  \int_{\Pi}|\nabla \phi_\Pi^i-j_i|^2~dx$. 
\[
\left| \int_{\Pi}|\nabla \phi_\Pi^i|^2~dx+ \int_{\Pi}j_i^2~dx-2\int_{\partial \mathcal{B}_i}(j_i\cdot n)~ds\right|=\int_{\Pi}|\nabla \phi_\Pi^i-j_i|^2~dx.
\]
This quantity is bounded, hence, the asymptotics of $\mathcal{G}_i$ is given by \eqref{E:en-test-flux}, whose asymptotics 
%\[
%%%%\int_{\Pi}|\nabla \phi_\Pi^i|^2~dx=[1+o(1)]\int_{\Pi}\frac{d\bar{x}}{H^2(\bar{x})}, \quad \mbox{when }~ \delta\ll 1,
%\int_{\Pi}j_i^2~dx= \int_{\Pi}\frac{d\bar{x}}{H^2(\bar{x})}, \quad \mbox{when }~ \delta\ll 1.
%\]
 in its turn is shown in Appendix \ref{A:coeff}, see  also \cite{bk,bn,bgn1}), and is given by:
\[
\int_{\Pi}\frac{d\bar{x}}{H^2(\bar{x})}=g_\delta[1+o(1)],\quad  \mbox{for } \delta\ll 1.
\]
$\Box$

\subsection{Constant $\mathcal{C}_{12}$ in definition of $g_\delta$} \label{A:coeff}

Here we show what is the constant $\mathcal{C}_{12}$ in asymptotics of $g_\delta$ that we claimed to be dependent on curvatures of particles boundaries at the point of the smallest distance between each other.

In the cylindrical coordinate system introduced above, that is, the one with the axis $Ox_d$ coinciding with the line of the closest distance between $\mathcal{B}_1$ and $\mathcal{B}_2$, and with the origin at the mid-point of this line, the boundaries $\partial\mathcal{B}_1$ and $\partial\mathcal{B}_2$ are approximated by parabolas ($d=2$) and paraboloids ($d=3$):
\begin{equation} \label{E:oscul}
\begin{array}{r l r l l}
\partial\mathcal{B}_1: & \displaystyle y=\frac{\delta}{2}+\frac{x^2}{2\alpha_1}, & \partial\mathcal{B}_2: & \displaystyle y=-\frac{\delta}{2}-\frac{x^2}{2 \alpha_2}, & d=2\\[8pt]
\partial\mathcal{B}_1: & \displaystyle z=\frac{\delta}{2}+\frac{x^2}{2a_1}+\frac{y^2}{2b_1}, & \partial\mathcal{B}_2: & \displaystyle z=-\frac{\delta}{2}-\frac{x^2}{2a_2}-\frac{y^2}{2b_2}, & d=3\\[5pt]
\end{array}
\end{equation}
The distance between these paraboloids is 
\begin{equation} \label{E:dist}
h(\bar{x})=\begin{cases}
\displaystyle  \delta+\frac{x^2}{\alpha}, \quad \alpha:=\frac{2\alpha_1\alpha_2}{\alpha_1+\alpha_2}, &  d=2\\[8pt]
\displaystyle  \delta+\frac{x^2}{a}+\frac{y^2}{b}, \quad  a:=\frac{2a_1a_2}{a_1+a_2},~ b:=\frac{2b_1b_2}{b_1+b_2}, & d=3
\end{cases}
%\begin{array}{r r l l l}
%h(\bar{x})=&h(x)& \displaystyle = \delta+\frac{x^2}{\alpha}, & \displaystyle  \alpha:=\frac{2\alpha_1\alpha_2}{\alpha_1+\alpha_2}, & d=2\\[8pt]
%h(\bar{x})=&h(x,y)& \displaystyle = \delta+\frac{x^2}{a}+\frac{y^2}{b}, & \displaystyle  a:=\frac{2a_1a_2}{a_1+a_2},~ b:=\frac{2b_1b_2}{b_1+b_2}, & d=3
%\end{array}
\end{equation}
For sufficiently small neck-width $w\ll 1$, this distance $h(\bar{x})$ by \eqref{E:dist} is a ``good'' approximation for the actual distance $H(\bar{x})$ between the boundaries $\partial\mathcal{B}_1$ and $\partial\mathcal{B}_2$ in the sense that
\[
\int_{\Pi}\frac{d\bar{x}}{H^2(\bar{x})}=\int_{\Pi}\frac{d\bar{x}}{h^2(\bar{x})}[1+o(1)], %\quad \mbox{for }~w\ll 1,
\]
that is, provides the leading asymptotics of $\mathcal{G}_i$ from Appendix \ref{A:G-lem}.
Going back to \eqref{E:dist}, we note that in 2D the parameter $\alpha$ is the harmonic mean of the radii of curvatures of parabolas approximating $\partial\mathcal{B}_1$ and $\partial\mathcal{B}_2$. Similarly, in 3D quantities $a$ and $b$ are related to the Gaussian $\mathcal{K}_i$ and mean $\mathcal{H}_i$ curvatures of the corresponding paraboloids at the points of the their closest distance via:
\[
\mathcal{K}_i=\frac{4}{a_ib_i}, \qquad \mathcal{H}_i=\frac{a_i+b_i}{a_ib_i}, \quad i=1,2.
\]
Finally, direct evaluating of the integral $\displaystyle \int_{\Pi}\frac{d\bar{x}}{h^2(\bar{x})}$ yields the main asymptotic term for $\mathcal{G}_i$ as $\delta \ll 1$ and defines $g_\delta$ of \eqref{E:g}:
\[
\int_{\Pi}\frac{d\bar{x}}{h^2(\bar{x})}=[1+o(1)]
\begin{cases}
\pi \alpha \delta^{1/2}, & d=2\\[3pt]
\pi a b |\ln \delta |, & d=3
\end{cases}%=g_\delta[1+o(1)]
 \quad \mbox{for }~ \delta\ll 1,
\]
where $\alpha$, $a$, $b$ are defined in \eqref{E:dist} in terms of coefficients of the osculating paraboloids \eqref{E:oscul} at the point of the closest distance between particles surfaces. 
Thus,
\begin{equation} \label{E:C0}
\mathcal{C}_{12}=
\begin{cases}
\pi \alpha , & d=2\\[3pt]
\pi a b , & d=3
\end{cases}
\end{equation}
$\Box$

\end{document}